# Geometry Mapping, Complete Pascal Scheme versus Standard Bilinear Approach

Sulaiman Y. Abo Diab,

*Department of Structural Mechanics, Tishreen University, Faculty of Civil Engineering*, Lattakia, Syria.
email:   sabodiab@tishreen.edu.sy

*This paper presents a complete Pascal interpolation scheme for use in the plane geometry mapping applied in association with numerical methods. The geometry of a domain element is approximated by a complete Pascal polynomial. The interpolation procedure is formulated in a natural coordinate system. It also presents the methodology of constructing shape functions of Pascal type and establishing a transformation relation between natural and Cartesian variables. The performance of the presented approach is investigated firstly by calculating the geometrical properties of an arbitrary quadrilateral cross-section like area and moments of inertia and comparing the results with the exact values and with those provided by the standard linear approach and a serendipity family approach. Secondly, the assessment of the scheme follows using a straight-sided, compatible quadrilateral finite element for plate bending of which geometry is approximated by a complete set of second order with six free parameters. Triangular and quadrilateral shaped plates with different boundary conditions are computed and compared with well-known results in the literature. The presented procedure is of general applicability for elements with curved edges and not limited to straight-sided edges in the framework of numerical methods.*

## 1   Introduction

The first step in a finite element approximation consists in approximating the geometry of the finite element. The most wide spread geometric discretization using quadrilateral finite elements amongst other elements falls roughly under main categories of function classes, like Pascal, Lagrangian and Serendipity basis function. The earliest finite elements like the linear triangular element, the standard bilinear element, and the ACM plate-bending element were formulated related to a Cartesian coordinate system. These elements are the first examples of using Pascal, Lagrange and Serendipity class, respectively and formulated by the pioneer developer of the finite element method (Argyris, 1954), (Turner, et al., 1956), (Clough and Tocher, 1965), (Ergatoudis et al., 1968), (Zienkiewicz and Taylor, 1989), and others. A derivation of the stiffness matrix of the rectangular plate-bending element were also carried out in (Melosh, 1961, 1963). Detailed information about the history of developing the first finite elements can be found in, (Clough, 1990), (Zienkiewicz, 1995), (Gupta, 1996), (Felippa, 2017). Nowadays, a huge number of finite elements with complex geometry and shape functions are developed. Unfortunately, only very few of them based on Pascal basis. The research papers (Arnold, 2011), (Rand, et al., 2014), (Natarajan and Shanka 2018), discuss the use of Serendipity family within a finite element approximation; in the later one, a survey about Serendipity finite elements is included. There are also many kinds of functions and polynomials used for constructing the shape functions.  A construction of spline spaces, known as NURBS (Non-uniform rational B-spline curves and surfaces), is used by many authors (Hughes, 2005), (Floater and Lai, 20016), (Li et al., 2016). A survey on the local refinable splines can be found in the last one. Trigonometric and complex power series (Milsted and Hutchinson, 1974), (Piltner, 2003), rational functions (Wachspress 1975), see for example (Sukamar, 2006), Legendre, Lobatto (Vu and Deek, 2006), Chebyshev polynomials (Chan and Warburton, 2015), are also widely used. Certainly, there is a huge number of publications point out at the complexity of the topic concerned. The cited papers mentioned here are only few examples of a very wide research field, which still attracts many researchers. A comparison of high-order curved finite elements is included in (Sevilla et al., 2011). Another of different higher order finite element schemes for the simulation of Lamb waves is included in (Willberg et al., 2012). For approximating the geometry of a quadrilateral element, there are different interpolation schemes used in the literature. Amongst them is to mention the widely used standard bilinear scheme (Ergatoudis et. al. 1968) and the Wachspress scheme (as rational function), see for example (Wachspress 1975) in (Sukamar, 2006), as well as the Serendipity class scheme (Malik and Bert, 2000). A survey about construction of shape functions on convex polygons can be found in (Sukamar, 2006).

The standard bilinear approach is used for mapping a physical domain (quadrilateral shaped) to a computational domain (bi-unit square). This scheme is adopted and studied in the most textbooks and research papers of engineering sciences, see for example (Washizu, 1982), (Reddy, 1984), (Weaver and Johnston 1987), (Pilkey and Wunderlich, 1994), (Cook et al., 2001), (Douglas et. al., 2002), (Huang, et. al., 2014), (Shu et. al.,2000). Two geometry-approximation schemes of a quadrilateral element are compared in this paper, the mentioned standard bilinear scheme (of Lagrange type) and the Pascal interpolation scheme published for the first time in (Abo Diab, 2017). This paper presents also a quadrilateral finite element based on a compatible approximation basis for the displacement and a complete basis for the geometry. A transformation relation between natural and Cartesian variables is established using a complete Pascal polynomial of second order with six free parameters. The displacement field is constructed such that the c1–continuity requirements between adjacent elements are exactly fulfilled. The construction of conform displacement field is presented in (Abo Diab, 2018), in details. Free vibration analysis of quadrilateral and triangular shaped thin plates is performed and compared with well-known results from the literature obtained using different computational methods (Leissa, 1969), (Andesrson, 1954 in: Stokey, 2002), (Gorman, 2012), (Dozio and Carrera, 2011), (Xing and Liu, 2009), (Xing et. al., 2010), (Shi, D. et. al. 2018). In the last cited paper, a comparison with FEM-solution is included. Although triangular plates are never the best examples for assessing the performance of a quadrilateral finite element, the current element seems to behave well even for different boundary conditions. In the following, the construction of a compatible approximation basis for the displacement and a complete basis for the geometry are shortly outlined since the detailed formulation is presented in two previously published open access papers, (Abo Diab, 2017), (Abo Diab, 2018), and are available to the reader.

## 2   Geometry approximation in the computational domain

In the following, whenever is not pointed out, Latin indices range over the Cartesian co-ordinates and indices between round brackets identify the nodal points. For example $^i$ ranges over $x^i$ $(i=1,2)$, where $^{(p)}$ denotes the number of the nodal points.

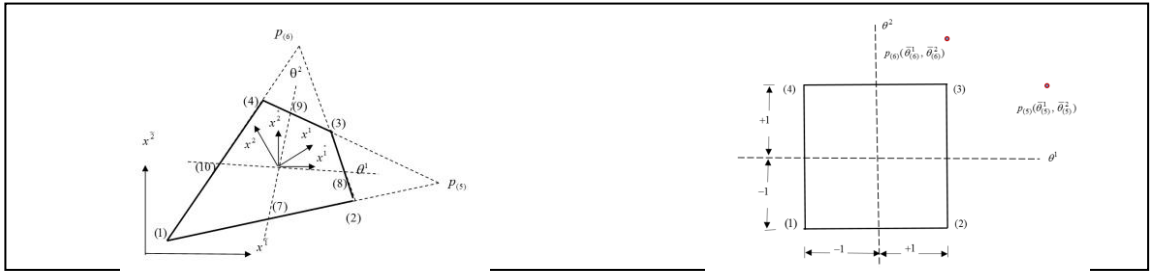

Figure 1.a. physical domain                Figure 1.b. Computational domain

Let $\Omega$ be a quadrilateral physical domain related to a Cartesian coordinate system $(x^{\tilde{1}}, x^{\tilde{2}})$ with the unit vectors $(e_{x^{\tilde{1}}}, e_{x^{\tilde{2}}})$ and defined by its four vertices $^{(1), (2), (3), (4)}$ (nodal points), see Figure 1.a.

$$x_{(p)}^{\tilde{i}} = \begin{bmatrix} x_{(1)}^{\tilde{1}} & x_{(1)}^{\tilde{2}} \\ x_{(2)}^{\tilde{1}} & x_{(2)}^{\tilde{2}} \\ x_{(3)}^{\tilde{1}} & x_{(3)}^{\tilde{2}} \\ x_{(4)}^{\tilde{1}} & x_{(4)}^{\tilde{2}} \end{bmatrix} \qquad (1)$$

Let the following values define the same nodal points in a natural coordinate system $(\theta^1, \theta^2)$. These nodal values span the same physical domain in what called the computational domain, Figure 1.b.

$$\theta^i_{(p)} = \begin{bmatrix} \theta^1_{(1)} & \theta^2_{(1)} \\ \theta^1_{(2)} & \theta^2_{(2)} \\ \theta^1_{(3)} & \theta^2_{(3)} \\ \theta^1_{(4)} & \theta^2_{(4)} \end{bmatrix} = \begin{bmatrix} -1 & -1 \\ +1 & -1 \\ +1 & +1 \\ -1 & +1 \end{bmatrix} \quad (2)$$

**2.1 Standard bilinear scheme**

The well-known standard bilinear interpolation approach involves polynomial with four undetermined parameters. Then, the scaled values Eqns. (2) of the natural coordinates at the four nodal points are sufficient to establish a transformation relation between the Cartesian variables and the natural variables.

The following transformation relates an arbitrary point with the natural coordinates $p'(\theta^1, \theta^2)$ to a corresponding point with the coordinates p ($x^{\tilde{1}}, x^{\tilde{2}}$)

$$\begin{bmatrix} x^{\tilde{1}} & x^{\tilde{2}} \end{bmatrix} = \frac{1}{4}\begin{bmatrix} (1-\theta^1)(1-\theta^2) & (1+\theta^1)(1-\theta^2) & (1+\theta^1)(1+\theta^2) & (1-\theta^1)(1+\theta^2) \end{bmatrix} * \begin{bmatrix} x^{\tilde{1}}_{(1)} & x^{\tilde{2}}_{(1)} \\ x^{\tilde{1}}_{(2)} & x^{\tilde{2}}_{(2)} \\ x^{\tilde{1}}_{(3)} & x^{\tilde{2}}_{(3)} \\ x^{\tilde{1}}_{(4)} & x^{\tilde{2}}_{(4)} \end{bmatrix} \quad (3)$$

Another form of the above transformation results in by defining the following differential geometry properties:

$$\begin{aligned}
a^{\tilde{1}}_{(1)} &= (x^{\tilde{1}}_{(1)} + x^{\tilde{1}}_{(2)} + x^{\tilde{1}}_{(3)} + x^{\tilde{1}}_{(4)})/4 & ; \quad a^{\tilde{2}}_{(1)} &= (x^{\tilde{2}}_{(1)} + x^{\tilde{2}}_{(2)}) + (x^{\tilde{2}}_{(3)} + x^{\tilde{2}}_{(4)})/4 \\
a^{\tilde{1}}_{(2)} &= ((x^{\tilde{1}}_{(2)} - x^{\tilde{1}}_{(1)}) - (x^{\tilde{1}}_{(4)} - x^{\tilde{1}}_{(3)}))/4 & ; \quad a^{\tilde{2}}_{(2)} &= ((x^{\tilde{2}}_{(2)} - x^{\tilde{2}}_{(1)}) - (x^{\tilde{2}}_{(4)} - x^{\tilde{2}}_{(3)}))/4 \\
a^{\tilde{1}}_{(3)} &= ((x^{\tilde{1}}_{(3)} - x^{\tilde{1}}_{(2)}) - (x^{\tilde{1}}_{(4)} - x^{\tilde{1}}_{(1)}))/4 & ; \quad a^{\tilde{2}}_{(3)} &= ((x^{\tilde{2}}_{(3)} - x^{\tilde{2}}_{(2)}) - (x^{\tilde{2}}_{(4)} - x^{\tilde{2}}_{(1)}))/4 \\
a^{\tilde{1}}_{(4)} &= ((x^{\tilde{1}}_{(3)} + x^{\tilde{1}}_{(1)}) - (x^{\tilde{1}}_{(4)} + x^{\tilde{1}}_{(2)}))/4 & ; \quad a^{\tilde{2}}_{(4)} &= ((x^{\tilde{2}}_{(3)} + x^{\tilde{2}}_{(1)}) - (x^{\tilde{2}}_{(4)} + x^{\tilde{2}}_{(2)}))/4
\end{aligned} \quad (4)$$

These parameters are called generalized parameters. The interpretation or the geometric meaning of the parameters $a^{\tilde{i}}_{(p)}$ is as follows:

$a^{\tilde{1}}_{(1)}, a^{\tilde{2}}_{(1)}$ are the Cartesian coordinates $(x^{\tilde{1}}, x^{\tilde{2}})$ of the geometric center of the element, respectively.

$a^{\tilde{1}}_{(2)}, a^{\tilde{2}}_{(2)}$ are a quarter of the difference of the projections of the edge [(1) (2)] and [(3) (4)] on the axes $(x^{\tilde{1}}, x^{\tilde{2}})$, respectively. They represent at the same time the $x^{\tilde{1}}$-components of the covariant base vectors evaluated at the geometric center of the element, which denoted as ($g_1^{x^{\tilde{1}}}$, $g_2^{x^{\tilde{1}}}$) in the definition (Eqn. (6)).

$a^{\tilde{1}}_{(3)}, a^{\tilde{2}}_{(3)}$ are a quarter of the difference of the projections of the edges [(2) (3)] and [(4) (1)] on the axes $(x^{\tilde{1}}, x^{\tilde{2}})$, respectively. They represent at the same time the $x^{\tilde{2}}$-components of the covariant base vectors evaluated at the geometric center of the element, which denoted as ($g_1^{x^{\tilde{2}}}$, $g_2^{x^{\tilde{2}}}$) in the definition (Eqn. (6))

$a^{\tilde{1}}_{(4)}, a^{\tilde{2}}_{(4)}$ are (in case of locating the origin of the Cartesian coordinate system at the element geometric center) a quarter of the difference between the projections of both diameters [(1) (3)] and [(4) (2)] on the axes $(x^{\tilde{1}}, x^{\tilde{2}})$, respectively. They represent at the same time the mixed derivatives of the covariant base vectors evaluated at the geometric center of the element, denoted as $g_{2,1}^{x^{\tilde{1}}}$, $g_{2,1}^{x^{\tilde{2}}}$ in Eqn. (7).

Thus, generalized parameters have the characters of deformation. The detailed formula of the transformation relation between natural and Cartesian variables using the generalized parameters reads:

$$\begin{aligned}
x^{\tilde{1}} &= a^{\tilde{1}}_{(1)} + a^{\tilde{1}}_{(2)}\theta^1 + a^{\tilde{1}}_{(3)}\theta^2 + a^{\tilde{1}}_{(4)}\theta^1\theta^2 \\
x^{\tilde{2}} &= a^{\tilde{2}}_{(1)} + a^{\tilde{2}}_{(2)}\theta^1 +, a^{\tilde{2}}_{(3)}\theta^2 + a^{\tilde{2}}_{(4)}\theta^1\theta^2
\end{aligned} \quad (5)$$

If we define the covariant base vectors and their derivatives as follows

$$\vec{g} = g_\alpha^{\tilde{i}} e_{\tilde{i}}; g_\alpha^{\tilde{i}} = \begin{bmatrix} g_1^{x^{\tilde{1}}} & g_1^{x^{\tilde{2}}} \\ g_2^{x^{\tilde{1}}} & g_2^{x^{\tilde{2}}} \end{bmatrix} = \begin{bmatrix} x^{\tilde{1}}_{,1} & x^{\tilde{2}}_{,1} \\ x^{\tilde{1}}_{,2} & x^{\tilde{2}}_{,2} \end{bmatrix} \quad (6)$$

$$\vec{g}_{\alpha,\beta} = g_{\alpha,\beta}^{\tilde{i}} e_{\tilde{i}}; g_\alpha^{\tilde{i}} = \begin{bmatrix} g_{1,1}^{x^{\tilde{1}}} & g_{1,1}^{x^{\tilde{2}}} \\ g_{2,1}^{x^{\tilde{1}}} & g_{2,1}^{x^{\tilde{2}}} \\ g_{1,2}^{x^{\tilde{1}}} & g_{1,2}^{x^{\tilde{2}}} \\ g_{2,2}^{x^{\tilde{1}}} & g_{2,2}^{x^{\tilde{2}}} \end{bmatrix} = \begin{bmatrix} x^{\tilde{1}}_{,11} & x^{\tilde{2}}_{,11} \\ x^{\tilde{1}}_{,21} & x^{\tilde{2}}_{,21} \\ x^{\tilde{1}}_{,12} & x^{\tilde{2}}_{,12} \\ x^{\tilde{1}}_{,22} & x^{\tilde{2}}_{,22} \end{bmatrix} \quad (7)$$

and evaluate the components of the above vectors at the geometric center, then we can express the transformation relation (Eqn. (5)) as follows:

$$
\begin{aligned}
x^{\tilde{1}} &= x^{\tilde{1}}_{(g)} + (g_1^{x^{\tilde{1}}})|_{(g)} \theta^1 + (g_2^{x^{\tilde{1}}})|_{(g)} \theta^2 + (g_{1,2}^{x^{\tilde{1}}})|_{(g)} \theta^1 \theta^2 \\
x^{\tilde{2}} &= x^{\tilde{2}}_{(g)} + (g_1^{x^{\tilde{2}}})|_{(g)} \theta^1 + (g_2^{x^{\tilde{2}}})|_{(g)} \theta^2 + (g_{1,2}^{x^{\tilde{2}}})|_{(g)} \theta^1 \theta^2
\end{aligned}
\quad (8)
$$

The notation $_{,1}$ denotes the derivative with respect to $\theta^1$ and the notation $_{,2}$ denotes the derivative with respect to $\theta^2$. The notation $()|_{(g)}$ means that the component is evaluated at the geometric center of the element.

Noting that $g_{1,2}^{\tilde{i}} = g_{2,1}^{\tilde{i}}$, then the last relation can be written in the following form:

$$
\begin{aligned}
x^{\tilde{1}} &= x^{\tilde{1}}_{(g)} + \begin{bmatrix} (g_1^{x^{\tilde{1}}})|_{(g)} & (g_2^{x^{\tilde{1}}})|_{(g)} \end{bmatrix} \begin{bmatrix} \theta^1 \\ \theta^2 \end{bmatrix} + \frac{1}{2}\begin{bmatrix} \theta^1 & \theta^2 \end{bmatrix} \begin{bmatrix} 0 & (g_{1,2}^{x^{\tilde{1}}})|_{(g)} \\ (g_{2,1}^{x^{\tilde{1}}})|_{(g)} & 0 \end{bmatrix} \begin{bmatrix} \theta^1 \\ \theta^2 \end{bmatrix} \\
x^{\tilde{2}} &= x^{\tilde{2}}_{(g)} + \begin{bmatrix} (g_1^{x^{\tilde{2}}})|_{(g)} & (g_2^{x^{\tilde{2}}})|_{(g)} \end{bmatrix} \begin{bmatrix} \theta^1 \\ \theta^2 \end{bmatrix} + \frac{1}{2}\begin{bmatrix} \theta^1 & \theta^2 \end{bmatrix} \begin{bmatrix} 0 & (g_{1,2}^{x^{\tilde{2}}})|_{(g)} \\ (g_{2,1}^{x^{\tilde{2}}})|_{(g)} & 0 \end{bmatrix} \begin{bmatrix} \theta^1 \\ \theta^2 \end{bmatrix}
\end{aligned}
\quad (9)
$$

This relation is similar to two variable Taylor expansion of the functions $x^{\tilde{i}}$ truncated after the second derivative terms, in the vicinity of the geometric center of the element, in which two non-linear terms are missed.

**2.2 Pascal interpolation scheme**

Consider the geometry mapping in a two-dimensional case. The current interpolation scheme makes use of Pascal polynomials. A complete Pascal polynomial of order or degree m contains all terms of the variables in that order and below. The order of a polynomial is defined by the highest power to which the variables are raised.

$$
\begin{array}{cccccc}
 & & 1 & & & \text{order } 0 \\
 & \theta^1 & & \theta^2 & & \text{order } 1 \\
 (\theta^1)^2 & & \theta^1\theta^2 & & (\theta^2)^2 & \text{order } 2 \\
 (\theta^1)^3 & (\theta^1)^2\theta^2 & & \theta^1(\theta^2)^2 & (\theta^2)^3 & \text{order } 3 \\
 (\theta^1)^4 & (\theta^1)^3\theta^2 & (\theta^1)^2(\theta^2)^2 & \theta^1(\theta^2)^3 & (\theta^2)^4 & \text{order } 4
\end{array}
\quad (10)
$$

The total number n of terms in a complete polynomial is given by the following relation:

$$n = \frac{1}{2}(m+1)(m+2) \quad (11)$$

For example, a fully quadratic mapping based on a complete Pascal interpolation scheme contains six terms with six undetermined parameters and a cubic mapping has ten terms. A quartic mapping contains fifteen terms and so on. At first, consider the quadratic mapping with six free parameters $a_{(p)}^{\tilde{i}}$ of the following form:

$$\begin{bmatrix} x^{\tilde{1}} & x^{\tilde{2}} \end{bmatrix} = \begin{bmatrix} 1 & \theta^1 & \theta^2 & (\theta^1)^2 & \theta^1\theta^2 & (\theta^2)^2 \end{bmatrix} \begin{bmatrix} a_{(1)}^{\tilde{1}} & a_{(1)}^{\tilde{2}} \\ a_{(2)}^{\tilde{1}} & a_{(2)}^{\tilde{2}} \\ a_{(3)}^{\tilde{1}} & a_{(3)}^{\tilde{2}} \\ a_{(4)}^{\tilde{1}} & a_{(4)}^{\tilde{2}} \\ a_{(5)}^{\tilde{1}} & a_{(5)}^{\tilde{2}} \\ a_{(6)}^{\tilde{1}} & a_{(6)}^{\tilde{2}} \end{bmatrix} \quad (12)$$

or in index form

$$x^{\tilde{i}} = M^{(p)} a_{(p)}^{\tilde{i}} \quad (13)$$

where $a_{(p)}^{\tilde{i}}$ are arbitrary parameters which have only a mathematical meaning at first. Compared with the standard bilinear approach, which contains only four terms $1$, $\theta^1$, $\theta^2$, $\theta^1\theta^2$, the quadratic terms $(\theta^1)^2$, $(\theta^2)^2$ are presented.

For determining the free parameters and giving them a geometrical meaning, a network of six nodes are needed in order to place a quadratic polynomial on the domain. For this purpose, the four element vertices (intersection point between every two neighboured edges) and the two poles, where every two opposite edges meet, are selected. The two poles belong to the geometry of the quadrilateral element as well. They play a key role in the approximation. In other words, more information than usually used in a standard bilinear approach must be passed to the interpolation procedure. The additional information include the Cartesian coordinates of poles $p_{(5)}$, $p_{(6)}$ and their natural coordinates. The Cartesian coordinates of the two intersections points of the lines crossing every two opposite edges (poles), i. e. $p_{(5)}$ intersection point of the two opposite edges $^{(1)\,(2)}$ and $^{(3)\,(4)}$ and $p_{(6)}$ intersection point of the other two opposite edges $^{(2)(3)}$ and $^{(4),\,(1)}$, can be computed easily by solving the system of the two equations representing the two edges, simultaneously.

Determining the natural coordinates of the two poles is not as easy as that of the Cartesian coordinates. A practical procedure for predicting the natural coordinates of the two poles $p_{(5)}$ and $p_{(6)}$ is described below:

The Cartesian coordinates of poles $p_{(5)}$ and $p_{(6)}$ are computed by solving the system of every two equations representing the opposite edges $^{(1)\,(2)}$, $^{(3)\,(4)}$ and $^{(2)(3)}$, $^{(4),\,(1)}$, simultaneously. As a result, we get the Cartesian coordinates $(x_{(5)}^{\tilde{1}}, x_{(5)}^{\tilde{2}}), (x_{(6)}^{\tilde{1}}, x_{(6)}^{\tilde{2}})$ of $p_{(5)}$ and $p_{(6)}$, respectively. The matrix $x_{(p)}^{\tilde{i}}$ will be extended to include these values:

$$x_{(p)}^{\tilde{i}} = \begin{bmatrix} x_{(1)}^{\tilde{1}} & x_{(1)}^{\tilde{2}} \\ x_{(2)}^{\tilde{1}} & x_{(2)}^{\tilde{2}} \\ x_{(3)}^{\tilde{1}} & x_{(3)}^{\tilde{2}} \\ x_{(4)}^{\tilde{1}} & x_{(4)}^{\tilde{2}} \\ x_{(5)}^{\tilde{1}} & x_{(5)}^{\tilde{2}} \\ x_{(6)}^{\tilde{1}} & x_{(6)}^{\tilde{2}} \end{bmatrix} \quad (14)$$

The matrix $\theta_{(p)}^{i}$ must be also extended to include the scaled values of the natural coordinates of the points $p_{(5)}$ and $p_{(6)}$. Let us denote these values as $(\theta_{(5)}^1, \theta_{(5)}^2), (\theta_{(6)}^1, \theta_{(6)}^2)$ and extend the natural nodal coordinates to include these values.

$$\theta^i_{(p)} = \begin{bmatrix} \theta^1_{(1)} & \theta^2_{(1)} \\ \theta^1_{(2)} & \theta^2_{(2)} \\ \theta^1_{(3)} & \theta^2_{(3)} \\ \theta^1_{(4)} & \theta^2_{(4)} \\ \theta^1_{(5)} & \theta^2_{(5)} \\ \theta^1_{(6)} & \theta^2_{(6)} \end{bmatrix} = \begin{bmatrix} -1 & -1 \\ +1 & -1 \\ +1 & +1 \\ -1 & +1 \\ \theta^1_{(5)} & \theta^2_{(5)} \\ \theta^1_{(6)} & \theta^2_{(6)} \end{bmatrix} \quad (15)$$

Now, the known Cartesian coordinates of the six nodes and the natural values of them with the four at first unknown values $(\theta^1_{(5)}, \theta^2_{(5)}), (\theta^1_{(6)}, \theta^2_{(6)})$ are passed to the interpolation procedure. Substituting the special natural coordinates of the nodal points of Eqn. (14) and (15) in Eqn. (12) yields, for details see Abo Diab (2017).

$$x^{\tilde{i}}_{(p)} = A^{(q)}_{(p)} a^{\tilde{i}}_{(q)} \quad (16)$$

where $A^{(q)}_{(p)}$ is a 6x6-matrix given by the following relation:

$$\begin{bmatrix} x^{\tilde{1}}_{(1)} & x^{\tilde{2}}_{(1)} \\ x^{\tilde{1}}_{(2)} & x^{\tilde{2}}_{(2)} \\ x^{\tilde{1}}_{(3)} & x^{\tilde{2}}_{(3)} \\ x^{\tilde{1}}_{(4)} & x^{\tilde{2}}_{(4)} \\ x^{\tilde{1}}_{(5)} & x^{\tilde{2}}_{(5)} \\ x^{\tilde{1}}_{(6)} & x^{\tilde{2}}_{(6)} \end{bmatrix} = \begin{bmatrix} 1 & -1 & -1 & 1 & 1 & 1 \\ 1 & 1 & -1 & 1 & -1 & 1 \\ 1 & 1 & 1 & 1 & 1 & 1 \\ 1 & -1 & 1 & 1 & -1 & 1 \\ 1 & \theta^1_{(5)} & \theta^2_{(5)} & (\theta^1_{(5)})^2 & \theta^1_{(5)}\theta^2_{(5)} & (\theta^2_{(5)})^2 \\ 1 & \theta^1_{(6)} & \theta^2_{(6)} & (\theta^1_{(6)})^2 & \theta^1_{(6)}\theta^2_{(6)} & (\theta^2_{(6)})^2 \end{bmatrix} \begin{bmatrix} a^{\tilde{1}}_{(1)} & a^{\tilde{2}}_{(1)} \\ a^{\tilde{1}}_{(2)} & a^{\tilde{2}}_{(2)} \\ a^{\tilde{1}}_{(3)} & a^{\tilde{2}}_{(3)} \\ a^{\tilde{1}}_{(4)} & a^{\tilde{2}}_{(4)} \\ a^{\tilde{1}}_{(5)} & a^{\tilde{2}}_{(5)} \\ a^{\tilde{1}}_{(6)} & a^{\tilde{2}}_{(6)} \end{bmatrix} \quad (17)$$

Solving the algebraic system of Eqns. (17) for the undetermined parameters gives the geometric meaning of them:

$$a^{\tilde{i}}_{(q)} = B^{(p)}_{(q)} x^{\tilde{i}}_{(p)} \quad (18)$$

where $B^{(p)}_{(q)}$ is the inverse matrix of $A^{(q)}_{(p)}$.

Substituting back the generalized parameters (18) into Eqns. (12), transformation relations between Cartesian and natural variables can be established:

$$x^{\tilde{i}} = M^{(q)} a^{\tilde{i}}_{(q)} = M^{(p)} B^{(q)}_{(p)} x^{\tilde{i}}_{(q)} = N^{(q)} x^{\tilde{i}}_{(q)} \quad (19)$$

where

$$N^{(q)} = M^{(p)} B^{(q)}_{(p)} \quad (20)$$

are the so called shape functions. It is of interest to note that the transformation relation before the geometric interpolation has the same form after the geometric interpolation except that the undetermined parameters $a^{\tilde{i}}_{(p)}$ are replaced by the generalized parameters $B^{(p)}_{(q)} x^{\tilde{i}}_{(p)}$. The parametric set $a^{\tilde{i}}_{(q)}$ is defined now by Eq. (18) and is dependent only on the constant known values of the Cartesian coordinates and first on the unknown values $(\theta^1_{(5)}, \theta^2_{(5)}), (\theta^1_{(6)}, \theta^2_{(6)})$.

The detailed formula of the transformation relation between natural and Cartesian coordinates using the generalized parameters is now given by the following relation:

$$\begin{aligned} x^{\tilde{1}} &= a^{\tilde{1}}_{(1)} + a^{\tilde{1}}_{(2)}\theta^1 + a^{\tilde{1}}_{(3)}\theta^2 + a^{\tilde{1}}_{(4)}(\theta^1)^2 + a^{\tilde{1}}_{(5)}\theta^1\theta^2 + a^{\tilde{1}}_{(6)}(\theta^2)^2 \\ x^{\tilde{2}} &= a^{\tilde{2}}_{(1)} + a^{\tilde{2}}_{(2)}\theta^1 +, a^{\tilde{2}}_{(3)}\theta^2 + a^{\tilde{2}}_{(4)}(\theta^1)^2 + a^{\tilde{2}}_{(5)}\theta^1\theta^2 + a^{\tilde{2}}_{(6)}(\theta^2)^2 \end{aligned} \quad (21)$$

Similarly, the transformation relation can also be expressed in terms of the components of the base vectors and their derivatives evaluated at the geometric center as follows:

$$x^{\tilde{1}} = x_{(g)}^{\tilde{1}} + (g_1^{x^{\tilde{1}}})|_{(g)} \theta^1 + (g_2^{x^{\tilde{1}}})|_{(g)} \theta^2 + (g_{1,1}^{x^{\tilde{1}}})|_{(g)} (\theta^1)^2 + (g_{1,2}^{x^{\tilde{1}}})|_{(g)} \theta^1 \theta^2 + (g_{2,2}^{x^{\tilde{1}}})|_{(g)} (\theta^2)^2$$
$$x^{\tilde{2}} = x_{(g)}^{\tilde{2}} + (g_1^{x^{\tilde{2}}})|_{(g)} \theta^1 + (g_2^{x^{\tilde{2}}})|_{(g)} \theta^2 + (g_{1,1}^{x^{\tilde{2}}})|_{(g)} (\theta^1)^2 + (g_{1,2}^{x^{\tilde{2}}})|_{(g)} \theta^1 \theta^2 + (g_{2,2}^{x^{\tilde{2}}})|_{(g)} (\theta^2)^2$$
(22)

Noting that $g_{1,2}^{\tilde{i}} = g_{2,1}^{\tilde{i}}$, then the last relation can be written in the following form:

$$x^{\tilde{1}} = x_{(g)}^{\tilde{1}} + \begin{bmatrix} (g_1^{x^{\tilde{1}}})|_{(g)} & (g_2^{x^{\tilde{1}}})|_{(g)} \end{bmatrix} \begin{bmatrix} \theta^1 \\ \theta^2 \end{bmatrix} + \frac{1}{2} \begin{bmatrix} \theta^1 & \theta^2 \end{bmatrix} \begin{bmatrix} (g_{1,1}^{x^{\tilde{1}}})|_{(g)} & (g_{1,2}^{x^{\tilde{1}}})|_{(g)} \\ (g_{2,1}^{x^{\tilde{1}}})|_{(g)} & (g_{2,2}^{x^{\tilde{1}}})|_{(g)} \end{bmatrix} \begin{bmatrix} \theta^1 \\ \theta^2 \end{bmatrix}$$

$$x^{\tilde{2}} = x_{(g)}^{\tilde{2}} + \begin{bmatrix} (g_1^{x^{\tilde{2}}})|_{(g)} & (g_2^{x^{\tilde{2}}})|_{(g)} \end{bmatrix} \begin{bmatrix} \theta^1 \\ \theta^2 \end{bmatrix} + \frac{1}{2} \begin{bmatrix} \theta^1 & \theta^2 \end{bmatrix} \begin{bmatrix} (g_{1,1}^{x^{\tilde{2}}})|_{(g)} & (g_{1,2}^{x^{\tilde{2}}})|_{(g)} \\ (g_{2,1}^{x^{\tilde{2}}})|_{(g)} & (g_{2,2}^{x^{\tilde{2}}})|_{(g)} \end{bmatrix} \begin{bmatrix} \theta^1 \\ \theta^2 \end{bmatrix}$$
(23)

This relation is similar to two variable Taylor expansion of the functions $x^{\tilde{i}}$ truncated after the second derivative terms, in the vicinity of the geometric center of the element, in which all non-linear second order terms are presented.

Eqns. (21) to (23) mean that the Cartesian variables $x^{\tilde{1}}, x^{\tilde{2}}$ in the physical domain are expressed as surfaces of the form $x^{\tilde{1}} = f(\theta^1, \theta^2), x^{\tilde{2}} = g(\theta^1, \theta^2)$ over the computational domain, where $f(\theta^1, \theta^2), g(\theta^1, \theta^2)$ 2D second order polynomial in $(\theta^1, \theta^2)$.

As may be seen, an explicit expression of the transformation relation between Cartesian and natural variables is first successful after determining the locations of the poles in the natural coordinate system.

It is observed during the inversion of $A_{(p)}^{(q)}$ that the generalized parameters $a_{(2)}^{\tilde{1}}, a_{(3)}^{\tilde{1}}, a_{(5)}^{\tilde{1}}, a_{(2)}^{\tilde{2}}, a_{(3)}^{\tilde{2}}, a_{(5)}^{\tilde{2}}$ remain constant and independent of the location of the poles in the natural coordinate system. This means in other words that the components of the covariant base vectors and the mixed second order derivatives evaluated at the geometric center retain their fixed values and independent of $(\theta_{(5)}^1, \theta_{(5)}^2), (\theta_{(6)}^1, \theta_{(6)}^2)$. The remaining generalized parameters $a_{(1)}^{\tilde{1}}, a_{(4)}^{\tilde{1}}, a_{(6)}^{\tilde{1}}, a_{(1)}^{\tilde{2}}, a_{(4)}^{\tilde{2}}, a_{(6)}^{\tilde{2}}$, expressed in terms of Cartesian coordinates of the geometric center and the components of the other two non-mixed second order derivatives of the covariant base vectors are functions of $(\theta_{(5)}^1, \theta_{(5)}^2), (\theta_{(6)}^1, \theta_{(6)}^2)$.

Exact determining of the natural coordinates of the pole $p_{(5)}$ can be curried out by substituting the expressions of $x^{\tilde{1}}, x^{\tilde{2}}$ from Eqn. (22) into the Cartesian edge equations of the edge (1) (2) and (3) (4) and solving the resulting equations with respect to $\theta^1, \theta^2$. Similarly, the natural coordinates of pole $p_{(6)}$ follow from solving the equations of the edges (2)(3) and (4), (1) expressed in terms of $\theta^1, \theta^2$. It is clear that the edge equations are quadratic in $\theta^1, \theta^2$ and therefore there exist multiple solutions. In addition, the numerical methods used in solving systems of two non-linear equations can be employed. Detailed discussion of solving non-linear equations are involved in many textbooks for Engineers and scientists, see for example (Hoffman, 1992).

An easy way to determine the pole coordinates $(\theta_{(5)}^1, \theta_{(5)}^2), (\theta_{(6)}^1, \theta_{(6)}^2)$ is to assume an initial guess corresponds to the zero curves of the tangent plane of the surfaces $x^{\tilde{1}} = f(\theta^1, \theta^2), x^{\tilde{2}} = g(\theta^1, \theta^2)$ and to improve those initial values iteratively (Newton's method). The location of the geometric center of the element offers investigating criteria for correct transformation relations.

The above discussion is valid for elements with curved edges, too. It is also possible to go higher in the approximation. For example, for a 2D third order Pascal interpolation scheme using ten term Pascal-polynomial, the four vertices and the two poles in addition to the four midpoints of the edges must be involved in the interpolation procedure. For a fourth order polynomial with fifteen free parameters, the parameters are determined from the natural coordinates of the four nodal points (vertices) of the element and the two intersections points of the lines crossing every two opposite edges as well as eight nodal points placed on the four edges and those of the geometric center of the element. The interpolation procedure leads to a 15x15 system of equations, which can be solved to find the undetermined parameters dependent on the Cartesian coordinates of the nodal points. But unfortunately the higher the approximation the higher is the computational cost in determining the natural coordinates of poles.

In the following, the procedure is investigated using a very simple example in calculating the area and moments of inertia of a quadrilateral cross section, for which exact values can be easily calculated.

## 3  Preminilary assessment

Example: A quadrilateral domain is defined by its four vertices [(1), (2) (3), (4)] (nodal points)

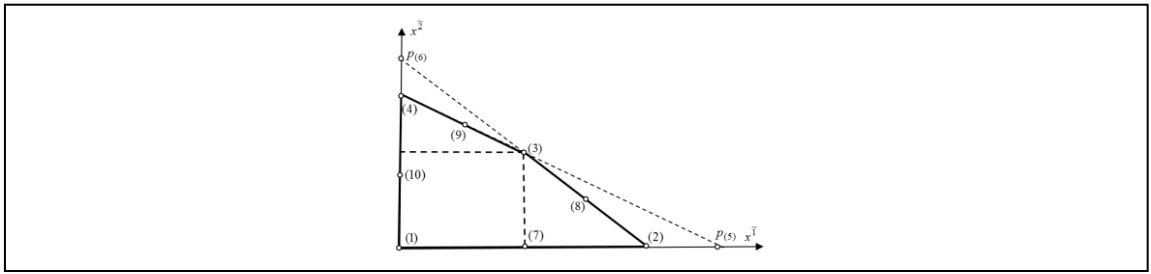

Figure 2 Quadrilateral cross-section, Cartesian coordinates of vertices, poles and edge midpoints

$$x^{\tilde{i}}_{(p)} = \begin{bmatrix} 0. & 0. \\ 8. & 0. \\ 4. & 3. \\ 0. & 5. \end{bmatrix} \quad (p) \downarrow$$

The exact values of the area of the quadrilateral cross section and the moments of inertia about the axes $(x^{\tilde{1}}, x^{\tilde{2}})$ can be computed easily as summation of these values of a rectangle with the side length 4x3 and two rectangled triangles with the catheti 4x3 and 4x2, Fig. 2 ( Timoshenko and Young, 1940)

$$A = (4)(3) + \frac{(4)(3)}{2} + \frac{(4)(2)}{2} = 22$$

$$I_{x^{\tilde{1}}} = \int_A (x^{\tilde{2}})^2 dA = \frac{4(3)^3}{12} + 4(3)(\frac{3}{2})^2 + \frac{4(3)^3}{36} + \frac{4(3)}{2}(\frac{1}{3}(3))^2 + \frac{4(2)^3}{36} + \frac{4(2)}{2}(3 + \frac{2}{3})^2 = 99.6667$$

$$I_{x^{\tilde{2}}} = \int_A (x^{\tilde{1}})^2 dA = \frac{3(4)^3}{12} + 4(3)(2)^2 + \frac{3(4)^3}{36} + \frac{4(3)}{2}(4 + \frac{4}{3})^2 + \frac{2(4)^3}{36} + \frac{2(4)}{2}(\frac{4}{3})^2 = 250.6667$$

$$I_{x^{\tilde{1}}x^{\tilde{2}}} = \int_A x^{\tilde{1}} x^{\tilde{2}} dA$$

$$= (4)(3)\frac{4}{2}\frac{3}{2} + (-\frac{(4)^2(3)^2}{72}) + \frac{(4)(3)}{2}(\frac{3}{3})(4 + \frac{4}{3}) + (-\frac{(2)^2(4)^2}{72}) + \frac{(2)(4)}{2}(3 + \frac{2}{3})(\frac{4}{3}) = 84.6667$$

**3.1 Standard bilinear approach**

Applying the standard bilinear approach yields the following transformation relation between the Cartesian coordinates and the natural element coordinates

$$\begin{bmatrix} x^{\tilde{1}} & x^{\tilde{2}} \end{bmatrix} = \frac{1}{4}\begin{bmatrix} N^{(1)} & N^{(2)} & N^{(3)} & N^{(4)} \end{bmatrix}\begin{bmatrix} x^{\tilde{1}}_{(1)} & x^{\tilde{2}}_{(1)} \\ x^{\tilde{1}}_{(2)} & x^{\tilde{2}}_{(2)} \\ x^{\tilde{1}}_{(3)} & x^{\tilde{2}}_{(3)} \\ x^{\tilde{1}}_{(4)} & x^{\tilde{2}}_{(4)} \end{bmatrix}$$

where

$$N^{(1)} = \frac{1}{4}(1-\theta^1)(1-\theta^2)$$

$$N^{(2)} = \frac{1}{4}(1+\theta^1)(1-\theta^2)$$

$$N^{(3)} = \frac{1}{4}(1+\theta^1)(1+\theta^2)$$

$$N^{(4)} = \frac{1}{4}(1-\theta^1)(1+\theta^2)$$

which reads after substituting the Cartesian coordinates of the element vertices

$$x^{\tilde{1}} = 3 + 3\theta^1 - \theta^2 - \theta^1\theta^2$$
$$x^{\tilde{2}} = 2 - 0{,}5\theta^1 + 2\theta^2 - 0.5\theta^1\theta^2$$

or in matrix form

$$x^{\tilde{1}} = 2 + \begin{bmatrix} 3 & -1 \end{bmatrix}\begin{bmatrix} \theta^1 \\ \theta^2 \end{bmatrix} + \frac{1}{2}\begin{bmatrix} \theta^1 & \theta^2 \end{bmatrix}\begin{bmatrix} 0 & -1 \\ -1 & 0 \end{bmatrix}\begin{bmatrix} \theta^1 \\ \theta^2 \end{bmatrix}$$

$$x^{\tilde{2}} = x^{\tilde{2}}_{(g)} + \begin{bmatrix} 2 & -0.5 \end{bmatrix}\begin{bmatrix} \theta^1 \\ \theta^2 \end{bmatrix} + \frac{1}{2}\begin{bmatrix} \theta^1 & \theta^2 \end{bmatrix}\begin{bmatrix} 0 & -0.5 \\ -0.5 & 0 \end{bmatrix}\begin{bmatrix} \theta^1 \\ \theta^2 \end{bmatrix}$$

Evaluating the integrals for calculating the element area and the moments of inertia in the computational domain using this transformation yields the exact values.

**3.2 Serendipity family approach**

Using a Serendipity family approximation in the form

$$\begin{bmatrix} x^{\tilde{1}} & x^{\tilde{2}} \end{bmatrix} = \begin{bmatrix} 1 & \theta^1 & \theta^2 & (\theta^1)^2 & \theta^1\theta^2 & (\theta^2)^2 & (\theta^1)^2\theta^2 & \theta^1(\theta^2)^2 \end{bmatrix}\begin{bmatrix} a^{\tilde{1}}_{(1)} & a^{\tilde{2}}_{(1)} \\ a^{\tilde{1}}_{(2)} & a^{\tilde{2}}_{(2)} \\ a^{\tilde{1}}_{(3)} & a^{\tilde{2}}_{(3)} \\ a^{\tilde{1}}_{(4)} & a^{\tilde{2}}_{(4)} \\ a^{\tilde{1}}_{(5)} & a^{\tilde{2}}_{(5)} \\ a^{\tilde{1}}_{(6)} & a^{\tilde{2}}_{(6)} \\ a^{\tilde{1}}_{(7)} & a^{\tilde{2}}_{(7)} \\ a^{\tilde{1}}_{(8)} & a^{\tilde{2}}_{(8)} \end{bmatrix}$$

and inserting the following Cartesian coordinates as well as the natural element coordinates of the edge-midpoints in addition to the vertices in the interpolation procedure

$$x_{(p)}^{\tilde{i}} = \begin{bmatrix} x_{(1)}^{\tilde{1}} & x_{(1)}^{\tilde{2}} \\ x_{(2)}^{\tilde{1}} & x_{(2)}^{\tilde{2}} \\ x_{(3)}^{\tilde{1}} & x_{(3)}^{\tilde{2}} \\ x_{(4)}^{\tilde{1}} & x_{(4)}^{\tilde{2}} \\ x_{(5)}^{\tilde{1}} & x_{(5)}^{\tilde{2}} \\ x_{(6)}^{\tilde{1}} & x_{(6)}^{\tilde{2}} \\ x_{(7)}^{\tilde{1}} & x_{(7)}^{\tilde{2}} \\ x_{(8)}^{\tilde{1}} & x_{(8)}^{\tilde{2}} \end{bmatrix} = \begin{bmatrix} 0.0 & 0.0 \\ 8.0 & 0.0 \\ 4.0 & 3.0 \\ 0.0 & 5.0 \\ 4.0 & 0.0 \\ 6.0 & 1.5 \\ 2.0 & 4.0 \\ 0.0 & 2.5 \end{bmatrix} ; \theta_{(p)}^{i} = \begin{bmatrix} -1.0 & -1.0 \\ 1.0 & -1.0 \\ 1.0 & 1.0 \\ -1.0 & 1.0 \\ 0.0 & -1.0 \\ 1.0 & 0.0 \\ 0.0 & 1.0 \\ -1.0 & 0.0 \end{bmatrix}$$

yield the transformation relation:

$$\begin{bmatrix} x^{\tilde{1}} & x^{\tilde{2}} \end{bmatrix} = \begin{bmatrix} N^{(1)} & N^{(3)} & N^{(3)} & N^{(4)} & N^{(5)} & N^{(6)} & N^{(7)} & N^{(8)} \end{bmatrix} \begin{bmatrix} x_{(1)}^{\tilde{1}} & x_{(1)}^{\tilde{2}} \\ x_{(2)}^{\tilde{1}} & x_{(2)}^{\tilde{2}} \\ x_{(3)}^{\tilde{1}} & x_{(3)}^{\tilde{2}} \\ x_{(4)}^{\tilde{1}} & x_{(4)}^{\tilde{2}} \\ x_{(5)}^{\tilde{1}} & x_{(5)}^{\tilde{2}} \\ x_{(6)}^{\tilde{1}} & x_{(6)}^{\tilde{2}} \\ x_{(7)}^{\tilde{1}} & x_{(7)}^{\tilde{2}} \\ x_{(8)}^{\tilde{1}} & x_{(8)}^{\tilde{2}} \end{bmatrix}$$

where

$$N^{(1)} = (-1 + (\theta^1)^2 + (\theta^2)^2 + \theta^1\theta^2 - (\theta^1)^2\theta^2 - \theta^1(\theta^2)^2)/4$$

$$N^{(2)} = (-1 + (\theta^1)^2 + (\theta^2)^2 - \theta^1\theta^2 - (\theta^1)^2\theta^2 + \theta^1(\theta^2)^2)/4$$

$$N^{(3)} = (-1 + (\theta^1)^2 + (\theta^2)^2 + \theta^1\theta^2 + (\theta^1)^2\theta^2 + \theta^1(\theta^2)^2)/4$$

$$N^{(4)} = (-1 + (\theta^1)^2 + (\theta^2)^2 - \theta^1\theta^2 + (\theta^1)^2\theta^2 - \theta^1(\theta^2)^2)/4$$

$$N^{(5)} = (1 - (\theta^1)^2 - \theta^2 + (\theta^1)^2\theta^2)/2$$

$$N^{(6)} = (1 + \theta^1 - (\theta^2)^2 - \theta^1(\theta^2)^2)/2$$

$$N^{(7)} = (1 + \theta^2 - (\theta^1)^2 - (\theta^1)^2\theta^2)/2$$

$$N^{(8)} = (1 - \theta^1 - (\theta^2)^2 + \theta^1(\theta^2)^2)/2$$

Finally, multiplying the shape function by the eight nodal coordinates $x_{(p)}^{\tilde{i}}$ gives the same transformation relation

$$x^{\tilde{1}} = 3 + 3\theta^1 - \theta^2 - \theta^1\theta^2$$
$$x^{\tilde{2}} = 2 - 0.5\theta^1 + 2\theta^2 - 0.5\theta^1\theta^2$$

Note that the higher order terms do not contribute to the transformation relation with any additional terms. They affect only the shape functions. Evaluating the integrals for element area and moments of inertia leads also to the exact values.

**3.2 Pascal approach**

Solving the equations of every two opposite edges of the quadrilateral shape, simultaneously give the Cartesian coordinates of their intersection points. Updating $x_{(p)}^{\tilde{i}}$ with the results gives the extended matrix of it.

$$x^{\tilde{i}}_{(p)} = \begin{bmatrix} 0.000 & 0.000 \\ 8.000 & 0.000 \\ 4.000 & 3.000 \\ 0.000 & 5.000 \\ 10.00 & 0.000 \\ 0.000 & 6.000 \end{bmatrix}$$

The following natural element coordinate sets of the four vertices and the two poles correspond to the Cartesian coordinates of the six nodes

$$\theta^i_{(p)} = \begin{bmatrix} \theta^1_{(1)} & \theta^2_{(1)} \\ \theta^1_{(2)} & \theta^2_{(2)} \\ \theta^1_{(3)} & \theta^2_{(3)} \\ \theta^1_{(4)} & \theta^2_{(4)} \\ \theta^1_{(5)} & \theta^2_{(5)} \\ \theta^1_{(6)} & \theta^2_{(6)} \end{bmatrix} = \begin{bmatrix} -1 & -1 \\ +1 & -1 \\ +1 & +1 \\ -1 & +1 \\ 4 & 1 \\ 1 & 3 \end{bmatrix} \quad \text{or} \quad \theta^i_{(p)} = \begin{bmatrix} \theta^1_{(1)} & \theta^2_{(1)} \\ \theta^1_{(2)} & \theta^2_{(2)} \\ \theta^1_{(3)} & \theta^2_{(3)} \\ \theta^1_{(4)} & \theta^2_{(4)} \\ \theta^1_{(5)} & \theta^2_{(5)} \\ \theta^1_{(6)} & \theta^2_{(6)} \end{bmatrix} = \begin{bmatrix} -1 & -1 \\ +1 & -1 \\ +1 & +1 \\ -1 & +1 \\ 3/2 & -1 \\ -1 & 7/5 \end{bmatrix}$$

Passing these values once a time for every set to the interpolation procedure described above gives the following shape functions:

$$\begin{bmatrix} x^{\tilde{1}} & x^{\tilde{2}} \end{bmatrix} = \begin{bmatrix} N^{(1)} & N^{(3)} & N^{(3)} & N^{(4)} & N^{(5)} & N^{(6)} \end{bmatrix} \begin{bmatrix} x^{\tilde{1}}_{(1)} & x^{\tilde{2}}_{(1)} \\ x^{\tilde{1}}_{(2)} & x^{\tilde{2}}_{(2)} \\ x^{\tilde{1}}_{(3)} & x^{\tilde{2}}_{(3)} \\ x^{\tilde{1}}_{(4)} & x^{\tilde{2}}_{(4)} \\ x^{\tilde{1}}_{(5)} & x^{\tilde{2}}_{(5)} \\ x^{\tilde{1}}_{(6)} & x^{\tilde{2}}_{(6)} \end{bmatrix}$$

The resulting shape functions for both natural element coordinate sets are listed below:

The shape functions produced by the first pole coordinate set are

$$N^{(1)} = (1 - \theta^1 - \theta^2 + \theta^1\theta^2)/4$$

$$N^{(2)} = (1/8 + \theta^1/4 - \theta^2/4 - \theta^1\theta^2/4 + (\theta^2)^2/4)$$

$$N^{(3)} = (2/3 + \theta^1/4 + \theta^2/4 + \theta^1\theta^2/4 - (\theta^1)^2/6 - (\theta^2)^2/4)$$

$$N^{(4)} = (3/20 - \theta^1/4 + \theta^2/4 - \theta^1\theta^2/4 + (\theta^1)^2/10)$$

$$N^{(5)} = (\theta^1)^2/15 - 1/15$$

$$N^{(6)} = (\theta^2)^2/8 - 1/8$$

The shape functions produced by the second pole coordinate set are:

$$N^{(1)} = (19/120 - \theta^1/4 - \theta^2/4 + \theta^1\theta^2/4 + (\theta^1)^2/15 + 5(\theta^2)^2/24)$$

$$N^{(2)} = (5/4 + \theta^1/4 - \theta^2/4 - \theta^1\theta^2/4 - (\theta^1)^2)$$

$$N^{(3)} = (1 + \theta^1 + \theta^2 + \theta^1\theta^2)/4$$

$$N^{(4)} = (3/2 - \theta^1/4 + \theta^2/4 - \theta^1\theta^2/4 - 5(\theta^2)^2/4)$$

$$N^{(5)} = 4(\theta^1)^2/5 - 4/5$$

$$N^{(6)} = 25(\theta^2)^2/24 - 25/24$$

Note that, substituting the natural coordinates of the four vertices and the two poles each once a time into the obtained transformation relation yields the Cartesian coordinates of these points. Substituting also the natural equation of an edge gives information about how this edge is bounded. The summation of the shape functions is equal to unity. It can be easily proved that the essential properties of shape functions are preserved.

Furthermore, the change in the approximation basis affects in the case of quadrilateral element with straight edges only the shape functions not the final transformation.

Multiplying the shape functions by the six nodal coordinates $\tilde{x}^{\tilde{i}}_{(p)}$ gives the same transformation relation as obtained previously using the standard bilinear approach and Serendipity family approach and consequently produces the same exact values for element area and moments of inertia...

## 4 Deflection and rotation approximation basis in the computational domain

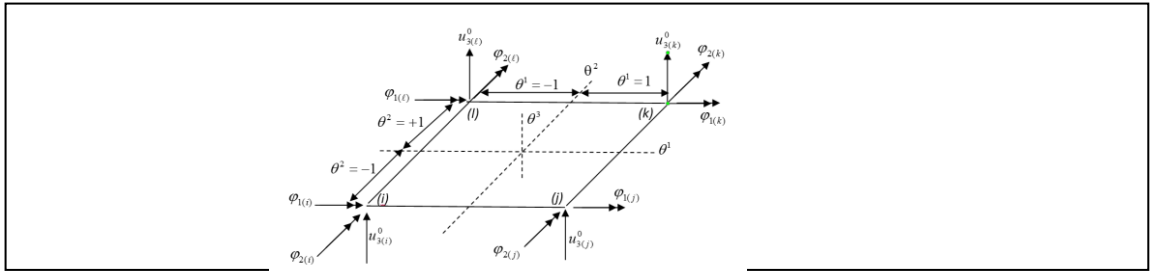

Figure 3 Bi-unit square, natural coordinate system, degrees of freedoms

The displacement approximation basis is constructed in the natural coordinate system. The element nodal displacement vector related to the natural coordinate system is as follow

$$u_{n(n)} = \{ u^0_{\theta^3(i)}, \varphi_{1(i)}, \varphi_{2(i)}, u^0_{\theta^3(j)}, \varphi_{1(j)}, \varphi_{2(j)}, u^0_{\theta^3(k)}, \varphi_{1(k)}, \varphi_{2(k)}, u^0_{\theta^3(\ell)}, \varphi_{1(\ell)}, \varphi_{2(\ell)} \} \qquad (24)$$

The displacements along the element boundary are approximated in terms of the classical Hermitian polynomials. From those displacements, the rotations are derived. Arranging the approximation functions for the rotations along the four element boundaries gives the following relation

$$\begin{bmatrix} \varphi_2^{(i)(j)} \\ \varphi_1^{(j)(k)} \\ \varphi_2^{(\ell)(k)} \\ \varphi_1^{(i)(\ell)} \end{bmatrix} = \begin{bmatrix} -h^1_{1,1} & 0 & -h^1_{2,1} & -h^1_{3,1} & 0 & -h^1_{4,1} & 0 & 0 & 0 & 0 & 0 & 0 \\ 0 & 0 & 0 & h^2_{1,2} & h^2_{2,2} & 0 & h^2_{3,2} & h^2_{4,2} & 0 & 0 & 0 & 0 \\ 0 & 0 & 0 & 0 & 0 & 0 & -h^1_{3,1} & 0 & -h^1_{4,1} & -h^1_{1,1} & 0 & -h^1_{2,1} \\ h^2_{1,2} & h^2_{2,2} & 0 & 0 & 0 & 0 & 0 & 0 & 0 & h^2_{3,2} & h^2_{4,2} & 0 \end{bmatrix} \begin{bmatrix} u^\circ_{\theta^3(i)} \\ \varphi_{1(i)} \\ \varphi_{2(i)} \\ u^\circ_{\theta^3(j)} \\ \varphi_{1(j)} \\ \varphi_{2(j)} \\ u^\circ_{\theta^3(k)} \\ \varphi_{1(k)} \\ \varphi_{2(k)} \\ u^\circ_{\theta^3(\ell)} \\ \varphi_{1(\ell)} \\ \varphi_{2(\ell)} \end{bmatrix} \qquad (25)$$

where

$$h_1^1 = \frac{1}{4}\left(2 - 3\theta^1 + \left(\theta^1\right)^3\right)$$
$$h_2^1 = \frac{1}{4}\left(-1 + \theta^1 + \left(\theta^1\right)^2 - \left(\theta^1\right)^3\right)$$
$$h_3^1 = \frac{1}{4}\left(2 + 3\theta^1 - \left(\theta^1\right)^3\right)$$
$$h_4^1 = \frac{1}{4}\left(1 + \theta^1 - \left(\theta^1\right)^2 - \left(\theta^1\right)^3\right)$$

(26)

and

$$h_1^2 = \frac{1}{4}\left(2 - 3\theta^2 + \left(\theta^2\right)^3\right);$$
$$h_2^2 = \frac{1}{4}\left(1 - \theta^2 - \left(\theta^2\right)^2 + \left(\theta^2\right)^3\right)$$
$$h_3^2 = \frac{1}{4}\left(2 + 3\theta^2 - \left(\theta^2\right)^3\right)$$
$$h_4^2 = \frac{1}{4}\left(-1 - \theta^2 + \left(\theta^2\right)^2 + \left(\theta^2\right)^3\right)$$

(27)

The notation $,_1$ denotes the derivative with respect to $\theta^1$ and the notation $,_2$ denotes the derivative with respect to $\theta^2$. The rotation $\varphi_1(\theta^1, \theta^2)$ inside the finite element will be interpolated depending on the rotation of every two opposite element boundaries $\varphi_1^{(j)(k)}(\theta^2)$, $\varphi_1^{(i)(\ell)}(\theta^2)$ as follows

$$\varphi_1(\theta^1, \theta^2) = \frac{1}{2}(1 - \theta^1)\varphi_1^{(i)(\ell)} + \frac{1}{2}(1 + \theta^1)\varphi_1^{(j)(k)}$$

(28)

Similarly, the rotation $\varphi_2(\theta^1, \theta^2)$ inside the finite element will be interpolated depending on the rotation of the other two opposite element boundaries $\varphi_2^{(i)(j)}(\theta^1)$, $\varphi_2^{(\ell)(k)}(\theta^1)$ as follows

$$\varphi_2(\theta^1, \theta^2) = \frac{1}{2}(1 - \theta^2)\varphi_2^{(i)(j)} + \frac{1}{2}(1 + \theta^2)\varphi_2^{(\ell)(k)}$$

(29)

Observing Eqn. (25) and Substituting into Eqns. (28) and (29) gives finally the following expression for the rotations

$$\begin{bmatrix}\varphi_1 \\ \varphi_2\end{bmatrix} = \begin{bmatrix} h_1^{01}*h_{1,2}^2 & h_1^{01}*h_{2,2}^2 & 0 & h_2^{01}*h_{1,2}^2 & h_2^{01}*h_{2,2}^2 & 0 & h_1^{01}*h_{3,2}^2 & h_1^{01}*h_{4,2}^2 & 0 & h_1^{01}*h_{2,2}^2 & h_1^{01}*h_{2,2}^2 & 0 \\ -h_1^{02}*h_{1,1}^1 & 0 & -h_1^{02}*h_{2,1}^1 & -h_1^{02}*h_{3,1}^1 & 0 & -h_1^{02}*h_{4,1}^1 & -h_2^{02}*h_{3,1}^1 & 0 & -h_2^{02}*h_{4,1}^1 & -h_2^{02}h_{1,1}^1 & 0 & -h_2^{02}*h_{2,1}^1 \end{bmatrix} \begin{bmatrix} u_{\theta^3(i)}^\circ \\ \varphi_{1(i)} \\ \varphi_{2(i)} \\ u_{\theta^3(j)}^\circ \\ \varphi_{1(j)} \\ \varphi_{2(j)} \\ u_{\theta^3(k)}^\circ \\ \varphi_{1(k)} \\ \varphi_{2(k)} \\ u_{\theta^3(\ell)}^\circ \\ \varphi_{1(\ell)} \\ \varphi_{2(\ell)} \end{bmatrix}$$

(30)

where

$$h_1^{01} = \frac{1}{2}(1 - \theta^1); \quad h_2^{01} = \frac{1}{2}(1 + \theta^1)$$
$$h_1^{02} = \frac{1}{2}(1 - \theta^2); \quad h_2^{02} = \frac{1}{2}(1 + \theta^2)$$

(31)

The approximation basis functions (31) are sufficient for developing the element stiffness matrix. The shape functions constructed in this way satisfy the $c^1$–continuity requirement between adjacent plate elements and the compatibility requirement with the classical beam and a relevant plane stress element when mixed at a folded plate edge (Müller and Abo Diab, 1987), Müller et al. (1987), Müller et al. (1991), (1994). But they show stiffened corner rotation. The corners become stiff and a fine mesh size is required in order to avoid this over stiffness. Another way consists in introducing the drilling degrees of freedom.

The explicit expression for the transverse displacement is necessary for developing the element load vector and the mass matrix. The following expansion of the deflection at any point $(\theta^1, \theta^2)$ is the simplest approximation of the deflection depending on the element rotations

$$u^0_{\theta^3}(\theta^1, \theta^2) = u^0_{\theta^3}(0,0) + \frac{\partial u^0_{\theta^3}}{\partial \theta^1}\theta^1 + \frac{\partial u^0_{\theta^3}}{\partial \theta^2}\theta^2 \qquad (32)$$

In Eqn. (32), $u^0_{\theta^3}(0,0)$ is the deflection at the geometric center of the element where $\theta^1 = 0, \theta^2 = 0$. When $\theta^1, \theta^2$ are small enough, the previous expression represents the linear terms of the multi variable Taylor expansion of the deflection of the point $(\theta^1, \theta^2)$ in the vicinity of the geometric center. It is also possible to go higher in the approximation by involving the higher order terms of the Taylor expansion. For example, considering the second term of the Taylor series gives

$$u^0_{\theta^3}(\theta^1, \theta^2) = u^0_{\theta^3}(0,0) + \frac{\partial u^0_{\theta^3}}{\partial \theta^1}\theta^1 + \frac{\partial u^0_{\theta^3}}{\partial \theta^2}\theta^2 + \frac{1}{2!}(\frac{\partial^2 u^0_{\theta^3}}{(\partial \theta^1)^2}(\theta^1)^2 + \frac{\partial}{\partial \theta^1}(\frac{\partial u^0_{\theta^3}}{\partial \theta^2})\theta^1\theta^2 + \frac{\partial}{\partial \theta^2}(\frac{\partial u^0_{\theta^3}}{\partial \theta^1})\theta^1\theta^2 + \frac{\partial^2 u^0_{\theta^3}}{(\partial \theta^2)^2}(\theta^2)^2) \qquad (33)$$

Assuming an average value for the deflection of the geometric center depending on deflections of the four nodal points,

$$\left[u^0_{\theta^3}(0,0)\right] = \begin{bmatrix} \frac{A_{(i)}}{A} & 0 & 0 & \frac{A_{(j)}}{A} & 0 & 0 & \frac{A_{(k)}}{A} & 0 & 0 & \frac{A_{(l)}}{A} & 0 & 0 \end{bmatrix} \begin{bmatrix} u^\circ_{\theta^3(i)} \\ \varphi_{1(i)} \\ \varphi_{2(i)} \\ u^\circ_{\theta^3(j)} \\ \varphi_{1(j)} \\ \varphi_{2(j)} \\ u^\circ_{\theta^3(k)} \\ \varphi_{1(k)} \\ \varphi_{2(k)} \\ u^\circ_{\theta^3(\ell)} \\ \varphi_{1(\ell)} \\ \varphi_{2(\ell)} \end{bmatrix} \qquad (34)$$

then we can evaluate the element load vector and the mass matrix in the usual way. In Eqn. (34), A is the element area and $A_{(p)}$; $(p)=(i),(j),(k),(l)$ is the subarea of the element included between the two element edges, which meet in (p), and the coordinate lines passing the geometric center. The selection of the deflection in the explained form is kinematically justified. By increasing the mesh refinement this selection represents a mathematically justified approximation for the deflection. Furthermore, in contrast to a Mindlin finite element application the deflection and the rotations are directly linked together. Note that the explicit formula of the deflection constructed corresponding to Eqns. (32), (33) is only used in order to predict the distribution of the load and mass over the element nodes. It does not affect in any way the stiffness matrix.

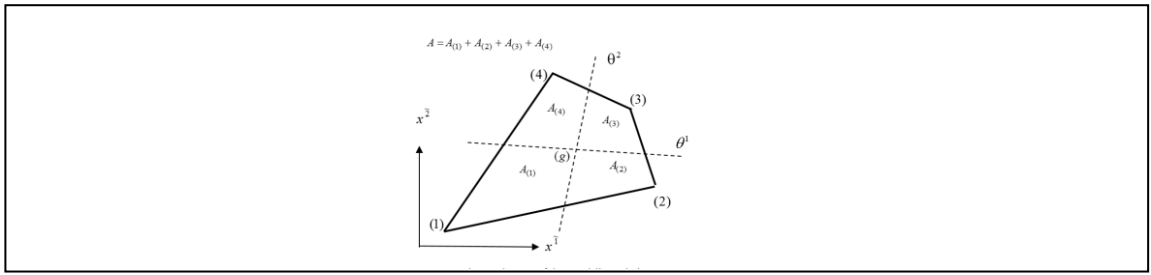

Figure 4: Subareas included between edges and coordinate lines

## 5  Variational approximation basis

The current finite element approximation is based on Hamilton's Principle. The 2D expression for the special case of the thin plate considered can be written in the absence of the prescribed boundary displacements in the following form:

$$\delta \int_{t_1}^{t_2} \left\{ \int_A \frac{1}{2} \chi_{ij} E^{ijkl} \chi_{kl} \, dA - \int_A \overline{q} \, u^0_{x^3} \, dA - \frac{1}{2} \int_A \dot{u}_i \rho^{ij} \dot{u}_j \, dA - \sum_{i=1}^{n} \overline{F}^{(i)} u^0_{x^3(i)} \right\} dt = 0 \tag{35}$$

where

$t_1$ and $t_2$ are two fixed time points of the vibration process, $\delta$ is the first variation, $E^{ijkl}$ is the tensor of the force-curvature dependency, $A$ is the element area and $dA$ its differential. $\dot{u}_i$ is the velocity vector in which both displacement and rotation components are included, $\rho^{ij}$ is the corresponding mass density matrix, $\overline{F}^{(i)}$ is the concentrated load applied at the point *(i)*. $\chi_{ij}$ is the curvature tensor. The derivation of the element is performed in the natural coordinate system, first with standard bilinear scheme and secondly with the Pascal scheme when approximating the geometry. The energy expression transformed into the natural coordinate system is as follow:

$$\delta \int_{t_1}^{t_2} \left\{ \int_A \frac{1}{2} \chi_{\alpha\beta} E^{\alpha\beta\gamma\delta} \chi_{\gamma\delta} \, dA - \int_A \overline{q} \, u^0_{\theta^3} \, dA - \frac{1}{2} \int_A \dot{u}_\alpha \rho^{\alpha\beta} \dot{u}_\beta \, dA - \sum_{i=1}^{n} \overline{F}^{(i)} u^0_{x^3(i)} \right\} dt = 0 \tag{36}$$

where

$$\chi_{\alpha\beta} = g^i_\alpha g^j_\beta \chi_{ij} \tag{37}$$

$$E^{\alpha\beta\gamma\delta} = g^\alpha_i g^\beta_j g^\gamma_k g^\delta_l E^{ijkl} \tag{38}$$

$g^i_\alpha$, $g^\alpha_i$ are the covariant and contra-variant components of the base vectors, respectively, (Meißner, 1986), (Klingbeil 1989)

The plate element is briefly described in (Abo Diab, 2018). Although the triangular plates are never the best examples for assessing a quadrilateral element, the numerical performance of the element will be shown by computing different triangular plates with different boundary conditions in addition to quadrilateral plates.

## 6  Numerical examples

The introduced conform thin plate bending element is verified numerically by convergence study of obtained solutions for triangular and quadrilateral plates with different boundary conditions. The geometry properties for all studied examples are given by the elastic modulus $E = 1365 \, \text{kN/m}^2$, the plate thickness t = 0.2 m, Poisson's ratio $\nu = 0.3$ and the mass density $\rho = 5 g / cm^3$ so that the flexural plate rigidity D and the factor $\sqrt{D/\rho h a^4}$ are kept equal to unity. The element stiffness and mass matrices are integrated numerically using a 3 × 3-point Gaussian rule. In all computed examples, the fundamental frequencies computed using the current element remain close to those listed from the literature but a significant difference in the higher frequencies can be observed. The solution obtained for the fundamental frequencies is slightly lower than that of the literature.

**6.1 Triangular plates**

**6.1.1 A cantilevered isosceles triangular plate**

The first computed triangular plate is a cantilevered isosceles one, clamped along the side $x^1 = 0$ (Fig. 5). The plate has the length to width ratio *b/a*=2, and the ratio *d/a*=0.5. The first four frequency parameters of the triangular plate $(\omega a^2)\sqrt{\rho h/D}$ for three different element mesh are listed in Tab. 1. The first line in the first cell is obtained by exact integration whereas the second line in the first cell is obtained using the numerical integration. For a comparison purpose, solutions obtained by (Leissa, 1969) and (Anderson, 1954) are also listed. The current solution is close to the listed solutions only for the fundamental frequency. It differs considerably for the higher frequencies. Some higher mode shapes show a plane-like vibration, which cannot be modeled by the beam theory. Figure 4a, 4b show the finite element mesh with 27 elements and the first six mode shapes, which corresponds to the first six frequencies, respectively.

Table 1 Comparison of the frequency parameter $(\omega a^2)\sqrt{\rho h/D}$ of a cantilevered isosceles triangular plate with the results provided in (Anderson,1954).

| mesh | $f_1$ | $f_2$ | $f_3$ | $f_4$ | $f_5$ | $f_6$ |
|---|---|---|---|---|---|---|
| 3 elements | *7.589261* | 26.034809 | 60.592412 | 73.075961 | 107.521903 | 250.155250 |
| 12 elements | 7.099645 | 32.794981 | 50.188296 | 73.257800 | 121.629890 | 133.980134 |
| 27 elements | 6.957378 | 30.474650 | 49.792043 | 74.532245 | 123.847674 | 134.889678 |
| (Leissa, 1969) | 7.122 | 30.718 | 90.105 | 259.4 | | |
| (Anderson, 1954) | 7.194 | 30.803 | 61.131 | 148.8 | | |

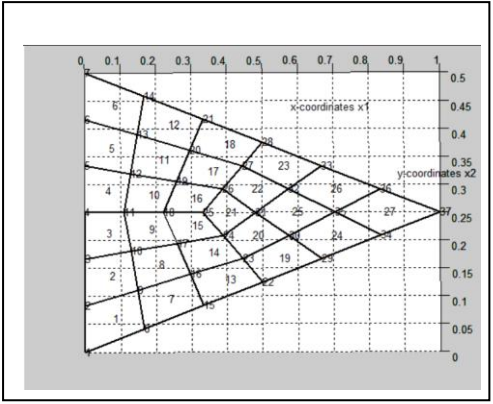

Figure 5a: A cantilevered isosceles triangular thin plate, geometry and finite element mesh.

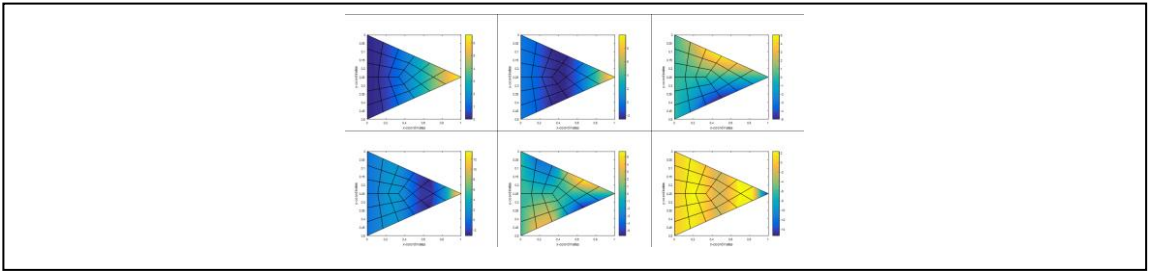

Figure 5b: The first six Eigen modes of a cantilevered isosceles triangular thin plate.

**6.1.2 A clamped isosceles triangular plate**

A clamped isosceles plate is subjected to modal analysis for a 3, 12 and 27 element mesh-size (Fig. 6). The mesh corresponds to the latter one is shown in figure 6a. The results are listed in Tab. 2 compared with a differential quadrature finite element solution provided by (Xing et. al., 2010). The results obtained using the current approach for the geometry are close to those provided by (Xing et. al., 2010) but they are still lower than those by about six to eight percent The less stiffener solution is obtained for the higher frequencies . A significant change in the higher frequencies can be observed using the Pascal scheme. Fig. 6b shows the mode shapes corresponds to frequencies listed in Tab.2.

Table 2: Comparison of the frequency parameter $(\omega a^2/\pi^2)\sqrt{\rho h/D}$ of isosceles triangular plate clamped along its all boundaries with the results provided in (Xing, et. al. 2010)

| mesh | $f_1$ | $f_2$ | $f_3$ | $f_4$ | $f_5$ | $f_6$ |
|---|---|---|---|---|---|---|
| 3 elements | 11.018331 | 14.876251 | 26.030531 | | | |
| 12 elements | 12.086224 | 13.817587 | 25.831078 | 27.882020 | 31.320584 | 36.864772 |
| 27 elements | 8.864298 | 14.947433 | 18.808521 | 22.624629 | 25.918113 | 27.944617 |
| Xing, (2010) | 9.502 | 15.987 | 19.735 | 24.600 | | |

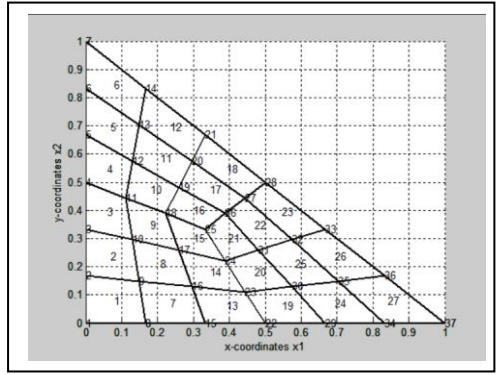

Figure 6a: All side clamped isosceles triangular plate meshed by 27-elements

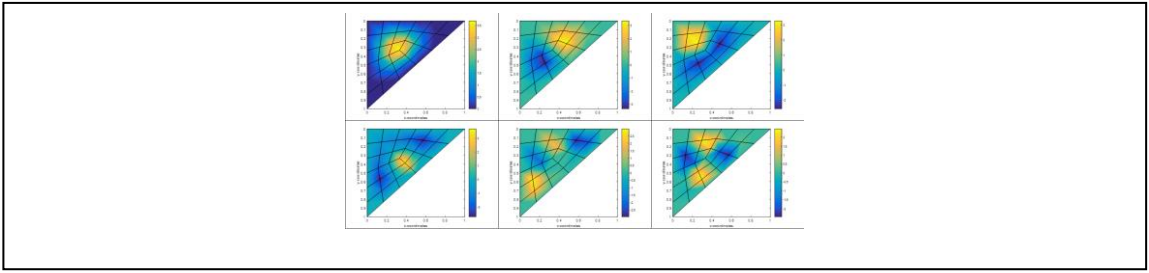

Figure 6b: The first six Eigen modes of a clamped isosceles triangular plate meshed by 27 elements

**6.1.3 A clamped equilateral triangular plate**

An equilateral triangular plate clamped along its all boundaries is subjected to modal analysis using various finite element meshes (Fig. 7). The results are listed in Tab. 3. A three-element mesh-size produces solution which is relatively close to that provided by (Xing, et. al., 2010), and (Leissa, 1969). The fundamental frequency is in a good agreement with the mentioned solutions. Figure 7b shows the mode shapes corresponds to the first six eigenvalues of the plate.

Table 3: Comparison of the frequency parameter $(\omega a^2/\pi^2)\sqrt{\rho h/D}$ of clamped equilateral triangular plate with the results provided in (Xing, et. al., 2010), (Leissa, 1969).

| mesh | $f_1$ | $f_2$ | $f_3$ | $f_4$ | $f_5$ | $f_6$ |
|---|---|---|---|---|---|---|
| 3 elements | 11.512539 | 17.226649 | 19.707365 | | | |
| 12 | 12.086224 | 13.817587 | 25.831078 | 27.882020 | 31.320584 | 36.864772 |
| elements | | | | | | |
| 27 elements | 9.541405 | 18.189182 | 18.497834 | 28.621356 | 29.341587 | 30.053141 |
| Leissa | $99.2/\pi^2$ | | | | | |
| (Xing and Liu 2010) | 9.502 | 15.987 | 19.735 | | | |

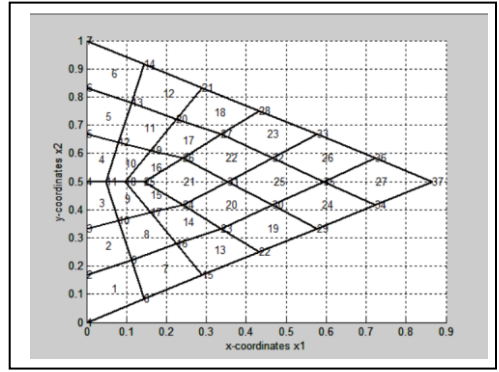

Figure 7a: An equilateral triangular plate clamped along its all boundaries meshed by 27 elements

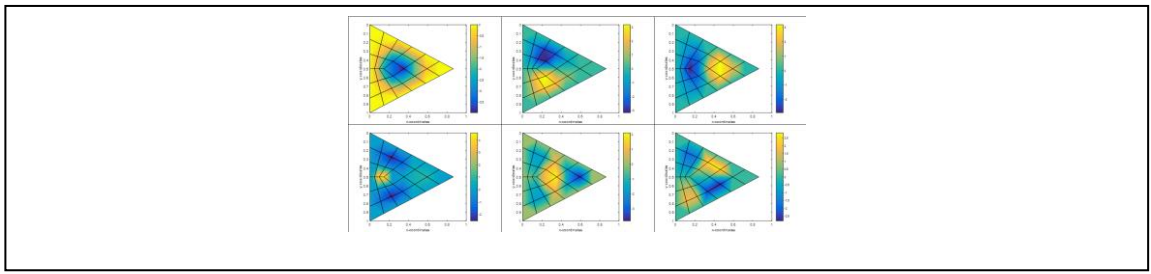

Figure 7b: The first six Eigen modes of a clamped equilateral triangular plate meshed by 27 elements

## 6.2 Quadrilateral plate

### 6.2.1. Clamped quadrilateral plate

Finally, the free vibration analysis of the fully clamped quadrilateral plate of fig. 8a were carried out using the present element. The four vertices of the plate are defined by the following global Cartesian coordinates:

$$\begin{bmatrix} x_{(1)}^{\tilde{1}} & x_{(1)}^{\tilde{2}} \\ x_{(2)}^{\tilde{1}} & x_{(2)}^{\tilde{2}} \\ x_{(3)}^{\tilde{1}} & x_{(3)}^{\tilde{2}} \\ x_{(4)}^{\tilde{1}} & x_{(4)}^{\tilde{2}} \end{bmatrix} = a \begin{bmatrix} 0. & 0. \\ 1. & 0. \\ 0.7929 & 0.7727. \\ 0.2394 & 0.6577 \end{bmatrix}$$

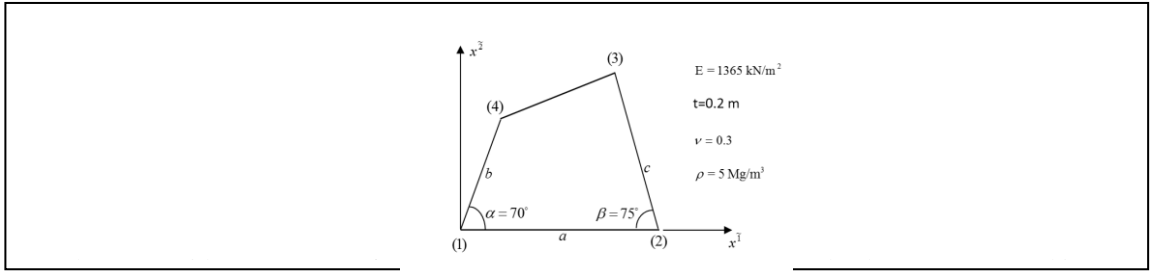

Figure 8a: Clamped quadrilateral plate, geometry, and material properties

The plate is meshed by different number of quadrilateral element. The results of a convergence study of the first 6.$^{\text{th}}$ non- dimensional natural frequencies $(\omega a^2/\pi^2)\sqrt{\rho h/D}$ using the proposed element are listed below in Tab.4. The present solution can be compared favorably with that provided by (Dozio and Carrera, 2011), for thickness-to-length ratios (h/a=0.005, 0.01, 0.05) and by (Shi, T. et. al. 2018).

Table 4: Comparison of the frequency parameter $(\omega a^2/\pi^2)\sqrt{\rho h/D}$ of fully clamped quadrilateral plate using the complete Pascal scheme

| mesh | $f_1$ | $f_2$ | $f_3$ | $f_4$ | $f_5$ | $f_6$ |
|---|---|---|---|---|---|---|
| 2x2 | 8.234961 | 11.329410 | 12.932349 | | | |
| 4x4 | 6.516769 | 12.481392 | 14.055341 | 18.910696 | 24.920696 | 27.287980 |
| 6x6 | 6.648953 | 12.623712 | 14.021059 | 19.526633 | 21.653798 | 25.592772 |
| 8x8 | 6.722913 | 12.759136 | 14.136358 | 19.856027 | 21.781022 | 25.444634 |
| (Dozio and Carrera, 2011) | 6.8294 | 13.115 | 14.319 | 20.777 | 22.296 | 25.627 |
| (Shi, D. et. al. 2018) | 6.816 | 13.096 | 14.300 | 20.753 | 22.262 | 25.597 |

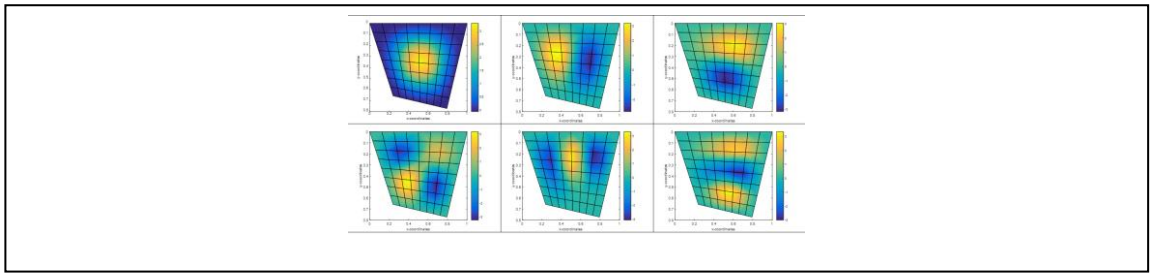

Figure 8b: The first six mode shapes of a clamped quadrilateral plate meshed by 64 elements

### 6.2.2. Quadrilateral cantilever plate

The last example is the general quadrilateral cantilever plate studied by (Dozio and Carrera, 2011) with side lengths a, b, c and internal angles $\alpha = 60°$, $\beta = 90°$ and the geometrical properties given in Fig. 9. The plate is subjected to modal analysis for a different mesh size. The four vertices of the plate are defined by the following global Cartesian coordinates:

$$\begin{bmatrix} x_{(1)}^{\tilde{1}} & x_{(1)}^{\tilde{2}} \\ x_{(2)}^{\tilde{1}} & x_{(2)}^{\tilde{2}} \\ x_{(3)}^{\tilde{1}} & x_{(3)}^{\tilde{2}} \\ x_{(4)}^{\tilde{1}} & x_{(4)}^{\tilde{2}} \end{bmatrix} = a \begin{bmatrix} 0 & 0 \\ 1. & 0 \\ 1. & 1. \\ 0.433 & 0.75 \end{bmatrix}$$

The first six frequency parameters of the cantilever quadrilateral plate $(\omega a^2/\pi^2)\sqrt{\rho h/D}$ for different mesh size are listed in Tab. 5. For the one element structure the element stiffness and mass matrices are integrated exactly, for other mesh size, the integration is performed numerically using a $3 \times 3$-point Gaussian rule. For a comparison purpose, a numerical solution obtained using a variable kinematic Ritz method applied to free vibration analysis of arbitrary quadrilateral thin and thick isotropic plates in (Dozio and Carrera, 2011) is also listed. Figure 9b shows the mode shapes corresponds to the first six eigenvalues of the plate whereas Figure 8c shows the same mode shapes in 3D-representation.

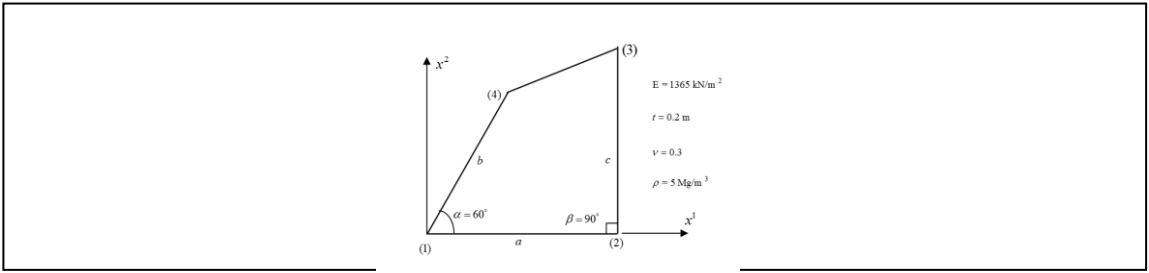

Figure 9a: quadrilateral cantilever thin plate, geometry and material properties

Table 5 Comparison of the frequency parameter $(\omega a^2/\pi^2)\sqrt{\rho h/D}$ of cantilever quadrilateral plate using the current element with the results provided in (Dozio and Carrera, 2011).

| mesh | $f_1$ | $f_2$ | $f_3$ | $f_4$ | $f_5$ | $f_6$ |
|---|---|---|---|---|---|---|
| 2x2 | 0.501838 | 1.596873 | 2.699458 | 4.134242 | 5.080273 | 6.210563 |
| 4x4 | 0.507972 | 1.641754 | 2.825560 | 4.369777 | 6.380308 | 7.684967 |
| 6x6 | 0.509355 | 1.655225 | 2.817862 | 4.385947 | 6.483664 | 7.676921 |
| 8x8 | 0.509880 | 1.660702 | 2.815623 | 4.399760 | 6.530852 | 7.670250 |
| Dozio and Carrera (2011) | 0.4857 | 1.3716 | 1.4692 | 2.2111 | 3.1421 | 3.3287 | 3.6373 | 4.5291 |

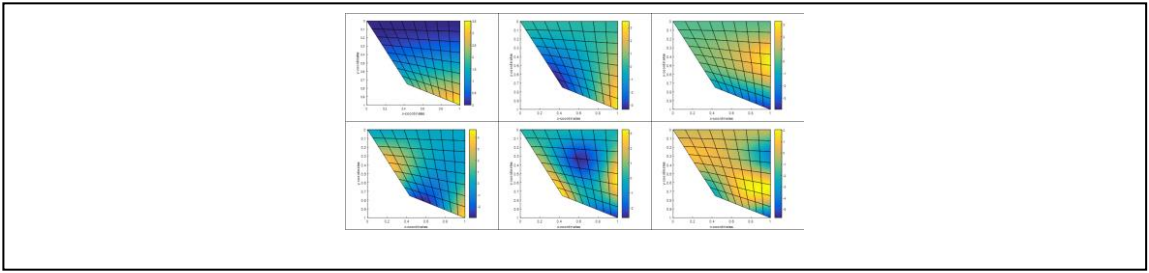

Figure 9b: The first six mode shapes of a cantilever quadrilateral plate meshed by 64 elements

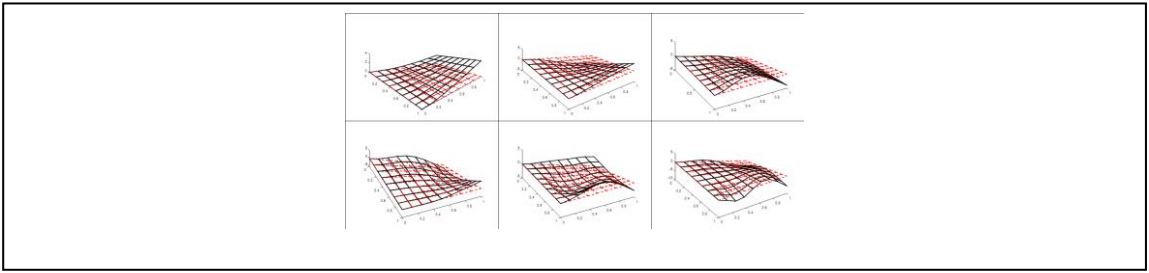

Figure 9c: 3D representation of the first six mode shapes of a cantilever quadrilateral plate meshed by 64 elements

# 7 Conclusion

A complete Pascal interpolation scheme for approximating the differential geometry properties of a quadrilateral shape and establishing a transformation relation between natural and Cartesian variables has been presented. The presented scheme is also compared with the standard bilinear approach and Serendipity family approach. It is found, that the different interpolation procedures are equivalent to a two variable Taylor expansion of the cartesian variables in a limited number of terms of the natural variables about the geometric center of the element, truncating the series after some non-linear terms. The different schemes led, in the case of a quadrilateral element, to different shape functions but finally to the same transformation relation. It is found that geometry approximation has also strong similarity to deformation problems in elasticity.

A compatible thin plate finite element with three degrees of freedoms at each node is presented. The element ensures continuity and inter-element continuity. Different quadrilateral and triangular plate with different boundary conditions are calculated. The results show that the element behaves well even for triangular plate.

_________________________________________________________________________________________

*Private Address*:   Sulaiman Abo Diab, Syria, Tartous, Hussain Al Baher
email:   sabodiab@hotmail.com